\documentclass[graybox,12pt]{svmult}

\usepackage{mathptmx}       
\usepackage{helvet}         
\usepackage{courier}        

\usepackage{makeidx}         
\usepackage{graphicx}        
\usepackage{multicol}        

\usepackage{amsmath}
\usepackage{amsfonts}
\usepackage{color}     
\usepackage{url}
\usepackage{verbatim}
\usepackage{enumerate}

\makeindex             

\usepackage{amssymb}
\usepackage{amsmath}
\usepackage{epsfig}
\usepackage{mathptmx}       
\usepackage{helvet}         
\usepackage{courier}        
\usepackage{type1cm}        

\usepackage{makeidx}         
\usepackage{graphicx}        
\usepackage{multicol}        
\usepackage[bottom]{footmisc}

\usepackage{amsbsy}
\usepackage{MnSymbol}
\usepackage{amsfonts}
\usepackage{amsopn}
\usepackage{dsfont}
\usepackage[english]{babel}
\usepackage[applemac]{inputenc}
\usepackage[colorlinks=true,breaklinks=true]{hyperref}
\usepackage[T1]{fontenc}
\usepackage{xcolor}
\usepackage{enumerate}
\usepackage{latexsym,amsmath,eufrak}
\usepackage{xcolor}
\usepackage{stackrel}
\newtheorem{teo}{Theorem}[section]
\newtheorem{prop}[teo]{Proposition}
\newtheorem{lem}[teo]{Lemma}
\newtheorem{defi}[teo]{Definition}
\newtheorem{rmk}[teo]{Remark}
\newtheorem{ex}[teo]{Example}

\newenvironment{proofp}{\noindent {\bf Proof.}}{\ \ \  \qed}

\newcommand\dstyle\displaystyle

\newcommand\R{\mathbb{R}}

\newcommand\al\alpha
\newcommand\be\beta
\newcommand\de\delta
\newcommand\la\lambda
\newcommand\tha\theta

\newcommand\iy\infty

\newcommand{\stk}{\stackrel{k}{\star}}

\newcommand{\hyp}[5]{\,\mbox{}_{#1}F_{#2}\!\left(
  \genfrac{}{}{0pt}{}{#3}{#4};#5\right)}

\newcommand\bma{\begin{pmatrix}}
\newcommand\ema{\end{pmatrix}}

\newmuskip\pFqskip
\mathchardef\pFcomma=\mathcode`, 

\begin{document}
\title*{On the Fractional Dunkl Laplacian}
\author{
Fethi Bouzeffour
\thanks{Universit\'e de Carthage, Facult\'e des Sciences de Bizerte, D\'epartement de Math\'ematiques,
Laboratoire d'Analyse Math\'ematiques et Applications LR11ES11.  2092 - El Manar I, Tunis, Tunisia.}
Wissem Jedidi
\thanks{Universit\'e de Tunis El Manar, Facult\'e des Sciences de Tunis, D\'epartement de Math\'ematiques,
Laboratoire d'Analyse Math\'ematiques et Applications LR11ES11.  2092 - El Manar I, Tunis, Tunisia}.
}
\maketitle
\abstract{In this paper, we present an approach to the fractional Dunkl Laplacian in a framework emerging from certain reflection symmetries in Euclidean spaces. Our main result is pointwise formulas, Bochner subordination, and an  extension problem  for the fractional Dunkl Laplacian as well.}
\section{Introduction}
In \cite{Dunkl}, Dunkl introduced a family of first-order differential-difference operators related to some finite reflection groups in the Euclidian space. Recently, these operators have gained considerable interest in various fields of mathematics and also in physical applications.  For more details about these operators see \cite{B1,B2,B3,Dunkl,11, KooBou,  Macdonald, Opdam, Roesler2,Roesler3, Trimeche} and references therein. The Dunkl-Laplacian operators are  $k$--deformations of the standard Laplacian operator $\Delta=\partial^2/\partial x_1^2+\,\dots\,+\partial^2/\partial x_d^2$ and they are fundamental tools for generalization of several  classical known results and so, it is a convenient setting for developing fractional Dunkl-Laplacian operators.  Recall that  for  a function $f$ in the space of  Schwartz functions $S(\R^d)$, the fractional Laplacian $(-\Delta)^{\alpha/2}, \;0<\alpha<2$, is defined  by means of the Fourier transform
$$(-\Delta)^{\alpha/2}f=\mathcal{F}^{-1}(|\xi|^{\alpha}\mathcal{F}f(\xi)),\quad  f\in S(\R^d),$$
and can be expressed by the pointwise formula \cite{Landkof}
\begin{equation}\label{R2}
(-\Delta)^{\alpha/2} f(x):= \frac{1}{\gamma_d(\alpha)}\lim_{\varepsilon\rightarrow 0}\,
\int_{\mathbb{R}^d\setminus B(0,\varepsilon)}\frac{f(x)-f(x-y)}{|y|^{d+\alpha}}\,dy,\quad
\mbox{where}\;\;\gamma_d(\alpha)= \frac{\pi^{d/2}|\Gamma(-\frac{\alpha}{2})|}{2^{\alpha}\Gamma(\frac{d+\alpha}{2})},
\end{equation}
and  $B(0, \varepsilon)$ denotes the ball of radius $\varepsilon$ centered at the origin, see \cite{Sting1}. The operator $(-\Delta)^{\alpha/2}$ is connected to PDE's through the Caffarelli--Silvestre extension theorem \cite{Caffarelli} which establishes the following:  if $U = U(x, y)$ is the solution to
\begin{equation}\label{ppb}
\left\{
  \begin{array}{l l}\displaystyle \Delta U(x,y)+\frac{\partial^2U(x,y)}{\partial^2y}+\frac{1-\alpha}{y}\frac{\partial U(x,y)}{\partial y}=0,\quad (x,y)\in \R^{d}\times (0,\,\iy),\\
  \\
  U(x,0)=u(x),\quad x\in\R^d,
  \end{array} \right.
       \end{equation}
then,
\begin{equation*}
(-\Delta)^{\alpha/2}u(x)=-\frac{2^{\alpha-1}
\Gamma(\alpha/2)}{\Gamma(1-\alpha/2)}\lim_{y\rightarrow 0^+}y^{1-\alpha}\frac{\partial U}{\partial y}(x,y).
\end{equation*}
The fractional Dunkl-Laplacian is then a natural object to consider, since it is a $k$-deformation of the standard Laplacian. Moreover, it is the simplest example of a large class of differential-difference operators associated with root systems. For more, we the refer to  Opdam's lecture notes  \cite{Opdam} for the trigonometric Dunkl theory, and to the books of Cherednik \cite{Cherednik} and of Macdonald \cite{Macdonald} for the generalized quantum theories.\\

In this work, we present several descriptions of the fractional Dunkl-Laplacian  on the Eucldian space. Our principal tools is the spherical mean-value type operator and Pizzetti's type formula related to the Dunkl operator. We are mainly interested in describing an analogue of  pointwise formula \eqref{R2} and the extension problem \eqref{ppb} in the setting of Dunkl theory.  In Section 2  we give a brief review of some elements of harmonic analysis related to the Dunkl operator. In Section 3 we derive several pointwise formulas  for the fractional Dunkl-Laplacian. In Section 4 we use the spherical mean-value type operator and Pezzitti's formula to give a Bochner type representation. In the last Section 5, we give  the fundamental solution and we study the fractional  Dirichlet-to-Neumann map for Dunkl-Laplacian.
\section{Preliminary}
In order to introduce our setting, we first collect some facts about the Dunkl operators. General references are \cite{Dunkl,11} and \cite{Roesler2,Roesler3}. Let $\mathcal{R}$ be a reduced root system in $\mathbb{R}^d$. For a vector $v \in \mathcal{R}$,  define the reflection $\sigma_v$ by
\begin{equation}
\sigma_v(x)=x-2\frac{\langle x,v\rangle}{|v|^2}v, \quad x\in \mathbb{R}^d,
\end{equation}
where $\langle\,.\,,\,.\,\rangle$ is the standard Euclidean inner product and $|x|=\sqrt{\langle x,x\rangle}$ is the Euclidean norm on $\R^d$. Let $k:\mathcal{R}\rightarrow[0,+\infty)$ be  a $G$-invariant function , where $G$ is the group of finite reflections related to root system $\mathcal{R}$ (the function $k$ is called multiplicity function). The Dunkl operators \,$T_j$,  $1\leq j\leq d$ are the following $k$--deformations of partial derivatives $\partial_j$ by difference operators:
\begin{align}\label{Dunkl}
T_j f(x)&=\partial_j f(x) +\frac{1}{2}\sum_{\,v\in\mathcal{ R}}\!k(v)\;\langle v,e_j\rangle\;
\frac{f(x)-f(\sigma_v(x))}{\langle v,\,x\rangle}\\
&=\partial_j f(x) +\sum_{\,v\in\mathcal{ R}_+}\!k(v)\;\langle v,e_j\rangle\; \frac{f(x)-f(\sigma_v(x))}{\langle v,\,x\rangle},\quad j=1,\,2,\dots,\,d. \nonumber
\end{align}
Here, $\mathcal{R}_+$ is any fixed positive subsystem of $\mathcal{R}$ and $e_1,\,\dots,\,e_d$, are the standard unit vectors of $\mathbb{R}^d.$  Note that the Dunkl operators $T_j$  commute pairwise and are skew-symmetric with respect to the $G$-invariant measure $w_k(x) dx$, where the  weight function  $w_k$ is given, for $x\in \R^d$,  by
\begin{equation*}
w_k(x):=\prod_{v\in \mathcal{R}} |\langle x,v\rangle|^{k(v)}=\prod_{v\in \mathcal{R}_+} |\langle x,v\rangle|^{2k(v)}.
\end{equation*}
The latter is a positive homogeneous function of degree $2\gamma_k$, where   $\gamma_k=\sum_{v\in\mathcal{R}_+}k(v)$. For fixed $x\in\R^d$, the Dunkl kernel $y\rightarrow\mathcal{E}_k(x,y)$ is the unique solution to the system
\begin{equation*}
\left\{
  \begin{array}{l l}\label{D1}
       T_j f=x_jf,\quad 1\leq j\leq d,\\
       \\
       f(0)=1.\end{array} \right.
\end{equation*}
The following integral formula was obtained by R\"{o}sler \cite{Roesler2}
\begin{equation}\label{EV}
\mathcal{E}_k(\xi,x)=\int_{\mathbb{R}^d}\hspace{-1mm}e^{\,\langle\xi,y\rangle}\,d\mu_x(y),
\quad x\!\in\!\mathbb{R}^d.
\end{equation}
where $\,\mu_x,\,x\!\in\!\mathbb{R}^d\hspace{.25mm},\;$  is a compactly supported probability measures. Specifically, it is supported in the convex hull $\mathcal{O}(x)$ of the $G$-orbit
of $x$. For a function $f$ in $L^1_k(\R^d)$, the Lebesgue space with respect to the measure $w_k(x)dx$,
the Dunkl transform  is defined by
\begin{equation}\label{Dunkltra}
\mathcal{F}_kf(\xi)=\widehat{f}(\xi)=\frac{1}{c_k} \int_{\mathbb{R}^d}\!f(x)\, \mathcal{E}_k(-i\,\xi,x)w_k(x)dx,\quad c_k\,=\int_{\mathbb{R}^d}\!e^{-\frac{|x|^2}2}\,w_k(x)dx.
\end{equation}
In a similar way to the Fourier transform (which is the particular case $k\equiv 0$), the Dunkl transform is a topological automorphism of the Schwartz space $S(\mathbb{R}^d)$ and can be extended to an isometric automorphism of $L^2_k(\R^d)$. Yet more, for every $f\!\in\!L^1_k(\R^d)$
such that $\mathcal{F}_kf\!\in\!L^1(\R^d, w_k)$, we have
$$f(\xi)=\mathcal{F}_k^2\!f(-\xi),\quad \xi\!\in\!\mathbb{R}^d,$$
and, for  $f\in S(\mathbb{R}^d)$,
\begin{equation}\label{dtj}
\mathcal{F}_k(T_jf)(\xi)=i\;\xi_j\;\mathcal{F}_kf(\xi),\quad \xi\!\in\!\mathbb{R}^d,\;\; 1\leq j\leq d.
\end{equation}
\par As in the classical case, a generalized translation operator is defined in the Dunkl setting side  on  $L^2(w_k)$ (the Lebesgue space of square integrable functions with respect to  $w_k(x)dx$) by Trim\`{e}che \cite{Trimeche}
\begin{eqnarray}\label{dutr}
\tau^xf(y):=\mathcal{F}_k^{-1}(\mathcal{E}_k(ix,y)\mathcal{F}_kf)(y), \quad y\in\mathbb{R}^d.
\end{eqnarray}
We also define the Dunkl convolution product for suitable functions $f$ and $g$ by
\begin{equation*}
f\stk g(x)=\frac{1}{c_k}\int_{\R^d}\tau^{-x}f(y)g(y)w_k(y)dy.
\end{equation*}
The Dunkl Laplacian associated with a reduced root system $\mathcal{R}$, and multiplicity function $k$, is the differential-difference operator, which acts on $C^2$ functions by
$$\Delta_k := \sum_{i=1}^ d T_{i}^2, \quad \mbox{\it where $T_i, \,1\leq i\leq d$ are the Dunkl operators defined in \eqref{Dunkl}. }$$
In explicit form, we have
$$\Delta_k f(x)=\Delta f(x)+2\sum_{v\in \mathcal{R}_+}k(v)\delta_vf(x),\quad f\in C^2(\mathbb{R}^d),$$
where $\Delta$ is the usual Laplacian on $\R^d$, and
\begin{equation*}
\delta_vf(x)=\frac{\langle\nabla f(x),v\rangle}{\langle v,x\rangle}-\frac{f(x)-f(\sigma_v(x))}{\langle v,x\rangle^2}.
\end{equation*}
Similarly to the fractional Laplacian on $\R^d$, the fractional powers of $(-\Delta_\kappa)^{\alpha/2}$
are defined by using the Dunkl transform \eqref{Dunkltra}. Indeed, the Dunkl Laplacian operator is essentially self-adjoint on $L^2_k(\R^d),$ see for instance \cite[Theorem 3.1]{Amri}. It is a Fourier-Dunkl multiplier with symbol $|\xi|^2,$ since by \eqref{dtj} we have
\begin{equation*}
\mathcal{F}_k(-\Delta_\kappa f)(\xi)=|\xi|^2 \mathcal{F}_k(f)(\xi).
\end{equation*}
Therefore, we can define in a natural way the action of $(-\Delta_k)^{\alpha/2}$  $(\alpha\in (0,\,2))$ on the Dunkl transform side by the equation
\begin{equation}\label{frac0}
\mathcal{F}_k((-\Delta_\kappa)^{\al/2} f)(\xi)=|\xi|^{\alpha} \mathcal{F}_k(f)(\xi),\quad \mbox{for all $f\in S(\R^d).$}
\end{equation}
\section{Fractional Dunkl Laplacian: Pointwise formulas}
\par Although we have formally introduced the fractional Dunkl-Laplacian by the formula \eqref{frac0}, such definition has a major disadvantage: it is not easy to understand a given function (or a distribution) by prescribing its Fourier Dunk transform. For this reason, we  introduce a different pointwise definition of the fractional Dunkl-Laplacian.
\begin{lem} For every $u\in S(\R^d),$ one has
\begin{equation*}
|2u(x)-\tau^x u(y)-\tau^{x} u(-y)|\leq  \frac{|y|^2}{c_k} \int_{\R^d}|\xi|^2 |\mathcal{F}_k u(\xi)|
\,w_k(\xi)\,d\xi.
\end{equation*}
\end{lem}
\begin{proofp} The integral representation for the Dunkl Kernel given by \eqref{EV}, leads to
$$|2-\mathcal{E}_k(-iy,\xi)-\mathcal{E}_k(-iy,\xi)|=2\int_{\R^d}\Big(1-\cos(\langle y,\eta\rangle)\Big)d\mu^k_\xi(\eta) \leq |y|^2\int_{\R^d}|\eta|^2d\mu^k_\xi(\eta).$$
Since the support of the probability measure $d\mu_\xi$ on $\R^d$ is contained in  the convex hull   $\mathcal{O}(\xi)$ of the orbit of $\xi$ under the action of the reflection group $G$, then
\begin{align*}
|2-\mathcal{E}_k(-iy,\xi)-\mathcal{E}_k(-iy,\xi)|\leq |y|^2|\xi|^2.
\end{align*}
On the other hand, the integral representation in \eqref{dutr} leads to
\begin{align*}
2u(x)-\tau^x u(y)&-\tau^{x} u(-y)= \frac{1}{c_k} \int_{\R^d} \mathcal{F}_k u(\xi)
\;\mathcal{E}_k(ix,\xi)\; (2-\mathcal{E}_k(iy, \xi)-\mathcal{E}_k(-iy, \xi)\; w_k(\xi) d\xi, \quad x,\, y\in \R^d.
\end{align*}
Therefore,
\begin{equation*}
|2u(x)-\tau^x u(y)-\tau^{x}u(-y)|\leq\frac{|y|^2}{c_k} \int_{\R^d}|\xi|^2|\mathcal{F}_k u(\xi)\,w_k(\xi)\,d\xi.
\end{equation*}
\end{proofp}

Recall that the Dunkl heat kernel $\Gamma_k(t,x)$ is given by  \cite{Roesler2}
$$\Gamma_k(t,x):=\, \frac{1}{(2t)^{\gamma_k +d/2}\;c_ k}\,e^{-|x|^2/4t},\quad x\in \mathbb{R}^d,\>t>0,$$
and has the  following properties:
\begin{equation} \label{f1}
\mathcal{F}_k\big(\Gamma_k(t,\,.\,)\big)(x)=e^{-t|x|^2}, \quad \int_{\mathbb{R}^d} \Gamma_k(t,x)\,w_k(y)=1,\quad t>0
\end{equation}
\begin{defi} Let $\alpha \in (0,2).$ The fractional Dunkl-Laplacian operator  $(-\Delta_k)^{\al/2}u$ of $u\in{S}(\R^d)$, is the nonlocal operator in $\R^d$ defined by
\begin{align}\label{fraa}
(-\Delta_k)^{\alpha/2} u(x)&=\,\frac{1}{\gamma_{k,d}(\alpha)}\int_{\mathbb{R}^d}
\frac{2u(x)-\tau^xu(y)-\tau^{-x}u(y)}{|y|^{\alpha+2\gamma_k+d}}w_k(y)\,dy.
\end{align}
where $\gamma_{k,d}(\alpha)$ is a suitable normalization constant that is given implicitly in Proposition \ref{4.3}.
\end{defi}
Notice that for $u\in S(\R^d),$ the integral in the right-hand side of \eqref{fraa} is convergent. Indeed, it suffices to write
\begin{eqnarray*}
\int_{\mathbb{R}^d} \Big|\frac{2u(x)-\tau^x u(y) -\tau^xu(-y)}{|y|^{\alpha+2\gamma_k+d}}\Big|w_k(y)dy
&=&\int_{\mathbb{R}^d \setminus B(0,\,\varepsilon)} \Big|\frac{2u(x)-\tau^x u(y) -\tau^xu(-y)}{|y|^{\alpha+2\gamma_k+d}}\Big|w_k(y)dy\\
&& \quad +\int_{B(0,\,\varepsilon)}\Big|\frac{2u(x)-\tau^x  u(y)-\tau^xu(-y)}{|y|^{\alpha+2\gamma_k+d}}\Big|w_k(y)dy.
\end{eqnarray*}
From Lemma 3.1, we have

$$\int_{ B(0,\,\varepsilon)}\Big|\frac{2u(x)-\tau^x u(y) -\tau^xu(-y)}{|y|^{\alpha+2\gamma_k+d}}
\Big|w_k(y)dy\leq C\int_{ B(0,\,\varepsilon)}\frac{w_k(y)}{|y|^{\alpha+2\gamma_k+d-2}} dy<\iy.$$
Since $w_k$ is a homogeneous function of degree $2\gamma_k,$ and $\alpha \in (0,2),$ the constant $C$ in the above inequality is given by $C=c^{-1}_k \big\|\,|\,.\,|^2\mathcal{F}_ku\,\big\|_1.$ On the other hand, keeping in mind that $u \in S (\R^d),$ which implies in particular that $\tau^xu \in  L^\iy_k(\R^d)$, then we have
\begin{align}
\int_{\mathbb{R}^d \setminus B(0,\,\varepsilon)}
\Big|\frac{2u(x)-\tau^x u(y) -\tau^xu(-y)}{|y|^{\alpha+2\gamma_k+d}}\Big|\;w_k(y)dy\leq 2\max(\|u\|_\iy,\,\|\tau^xu\|_\iy)\int_{\mathbb{R}^d \setminus B(0,\,\varepsilon)}
\frac{w_k(y)}{|y|^{\alpha+2\gamma_k+d}}dy<\iy.
\end{align}
Therefore Definition 3.3 provides a well-defined function on $\R^d$.

\begin{rmk} The following alternative expression for $(-\Delta_k)^{\al/2}$  is quite useful in the computations: For every $u\in S(\R^d)$, one has
\begin{equation}\label{Frr2}
(-\Delta_k)^{\alpha/2} u(x):=\,\frac{2}{\gamma_{k,d}(\alpha)}\text{PV}
\int_{\mathbb{R}^d} \frac{u(x)-\tau^xu(y)}{|y|^{\alpha+2\gamma_k+d}}\,w_k(y)dy,
\end{equation}
where $\text{PV}\int_{\mathbb{R}^d} v(y) dy=\lim_{\varepsilon\rightarrow 0} \int_{\mathbb{R}^d\setminus B(0,\varepsilon)}
 v(y) dy$ (Cauchy's principal value sense). Indeed,
\begin{align*}
\int_{\mathbb{R}^d}
\frac{2u(x)-\tau^x u(y) -\tau^xu(-y)}{|y|^{\alpha+2\gamma_k+d}}w_k(y)dy&=\lim_{\varepsilon\rightarrow 0^+} \int_{\mathbb{R}^d\setminus B(0,\,\varepsilon)} \frac{2u(x)-\tau^x u(y) -\tau^xu(-y)}{|y|^{\alpha+2\gamma_k+d}}w_k(y)dy\\
&= \lim_{\varepsilon\rightarrow 0^+} \int_{\mathbb{R}^d\setminus B(0,\,\varepsilon)}
\frac{u(x)-\tau^x u(y) }{|y|^{\alpha+2\gamma_k+d}}w_k(y)dy + \lim_{\varepsilon\rightarrow 0^+} \int_{\mathbb{R}^d\setminus B(0,\,\varepsilon)}
\frac{u(x)-\tau^x u(-y)}{|y|^{\alpha+2\gamma_k+d}}w_k(y)dy\\
&= 2\lim_{\varepsilon\rightarrow 0^+} \int_{\mathbb{R}^d\setminus B(0,\,\varepsilon)}
\frac{u(x)-\tau^x u(y)}{|y|^{\alpha+2\gamma_k+d}}w_k(y)dy.
\end{align*}
Note that it is  necessary to take the principal value of the last integral since we
have eliminated the cancellation of the linear terms in the symmetric difference of
order two, and $|u(x) -\tau^xu(y)|$ is only $O(|y|)$ (see \cite[Theorem 3.14]{ThangaveluXu}). Thus, the smoothness of the function $u$ no longer guarantees the local integrability.
\end{rmk}
\section{Representation via spherical mean-value operator}
In this section, we provide a useful expression of $(-\Delta_k)^{\alpha/2}u$ in terms of an integral involving the spherical mean-value operator associated to Dunkl operator. According  to \cite{Mejjaoli}, the spherical mean-value operator associated to Dunkl operator $\mathcal{M}_{r}^ku(x)$ is defined by
\begin{equation*}
\mathcal{M}_{r}^ku(x)=\frac{1}{\sigma_k(d)}
\int_{\mathbb{S}^{d-1}}\tau^xu(r\omega)w_k(\omega)d\sigma(\omega),\quad x\in \mathbb{R}^d,\,\,r\geq 0,
\end{equation*}
where $\mathbb{S}^{d-1}$ is the unit sphere in $\R^d$ , $d\sigma$ denotes the Lebesgue surface measure and
\begin{equation*}
\sigma_k(d):=\int_{\mathbb{S}^{d-1}}w_k(\omega)d\sigma(\omega) =\frac{c_k}{2^{\gamma_k+\frac{d}{2}-1} \Gamma(\gamma_k+\frac{d}{2})}.
\end{equation*}
From Mejjaoli-Trimeche \cite{Mejjaoli} one can extract the following proposition, that gives an extended Pizzetti's formula associated with the Dunkl operators.
\begin{prop}\label{mejme}
Let $f\in C^{\iy}(\R^d) $ and $a\in\R^d.$ The following asymptotic expansion is valid
\begin{equation}
\frac{1}{\sigma_k(d)}\int_{\mathbb{S}^{d-1}} \tau^af(\varepsilon\omega)\,w_k(\omega)d\sigma(\omega) \sim\Gamma(\gamma_k+\frac{d}{2}) \sum_{n=0}^{\iy}(\frac{\varepsilon}{2})^{2n}\frac{\Delta_k^nf(a)}
{n!\,\Gamma(\gamma_k+\frac{d}{2}+n)},\quad \mbox{as}\;\varepsilon\rightarrow 0^+.
\end{equation}
\end{prop}
\begin{lem}Let $u\in S(\R^d)$ one has
\begin{equation}\label{Intf}
(-\Delta_k)^{\al/2}u(x)=\frac{1}{\pi_{k,d}(\alpha)} \int_{0}^{\iy}\frac{u(x)-\mathcal{M}^k_ru(x)}{r^{1+\al}}\,dr,\quad x\in \R^d,\quad \mbox{\it where}\;\;\pi_{k,d}(\al):=\frac{ \gamma_{k,d}(\al)}{2\sigma_{k}(d)}.
\end{equation}
\end{lem}
\begin{proofp} From  Proposition \ref{mejme} (see also \cite[Theorem 4.17 ]{Mejjaoli}, one has
\begin{equation*}
\Delta_ku(x)=\frac{2(2\gamma_k+d)}{\sigma_k(\alpha)}\lim_{r\rightarrow 0}\frac{u(x)-\mathcal{M}_r^ku(x)}{r^2}.
\end{equation*}
Therefore, the integrand in the right-hand side of \eqref{Intf} behaves like
\begin{equation*}
\frac{u(x)-\mathcal{M}_r^ku(x)}{r^{1+r}}=O(r^{1-\alpha}),\quad r\rightarrow 0.
\end{equation*}
Since $u\in S(\R^d)$ and $\alpha \in (0,2),$ we conclude that the integral in the right-hand side of \eqref{Intf} is convergent. On the other hand, the use of the polar coordinates $x=r\omega,$ allows us to rewrite  \eqref{Frr2}, in the following forms
\begin{align*}
(-\Delta_k)^{\alpha/2} u(x)&=\frac{2}{\gamma_{k,d}(\alpha)} \lim_{\varepsilon\rightarrow 0}\,
\int_\varepsilon^\iy\frac{1}{r^{1+\alpha}} \int_{\mathbb{S}^{d-1}}(u(x)-\tau^xu(r\omega)) w_k(\omega)d\sigma (\omega)\,dr\\
&=\frac{2\sigma_{k}(d)}{ \gamma_{k,d}(\al)}\lim_{\varepsilon\rightarrow 0^+}\int_{\varepsilon}^{\iy} \frac{u(x)-\mathcal{M}_{r}^ku(x)}{r^{1+\al}}\,dr=\frac{1}{\pi_{k,d}(\alpha)} \int_{0}^{\iy}\frac{u(x)-\mathcal{M}^k_ru(x)}{r^{1+\al}}\,dr.
\end{align*}
 \end{proofp}
\subsection{Computation of the constant $\gamma_{k,d}(\alpha)$}
The reason behind the introduction of  the constant $\gamma_{k,d}(\alpha)$ in \eqref{fraa}, is to insure the validity of identity \eqref{frac0}. Its exact form is given here.
\begin{prop} The constant $\gamma_{k,d}(\alpha)$ in \eqref{fraa} is  given by
\begin{equation}\label{Cont}
\gamma_{k,d}(\alpha):=\frac{c_k|\Gamma(-\frac{\al}{2})|}{2^{\alpha +\gamma_k+d/2}
\Gamma(\gamma_k+\tfrac{\alpha+d}{2})}.
\end{equation}
Then, for every $u\in S(\R^d)$, we have
\begin{equation}\label{Fourf}
\mathcal{F}_k (-\Delta_k)^{\alpha/2}u(\xi)=|\xi|^\alpha \mathcal{F}_ku(\xi),\quad \xi \in \R^d.
\end{equation}
\label{4.3}\end{prop}
\begin{proofp} From \cite[formula 4.4]{Roesler4}, we have \begin{equation*}
\mathcal{M}_{r}^ku(x)=\frac{1}{c_k}\int_{\R^d}\mathcal{J}_{\gamma_k+d/2-1}(r|\xi|) \; \mathcal{F}_ku(\xi)\;\mathcal{E}_k(i\xi,x)\; w_k(\xi)\;d\xi,
\end{equation*}
where, the normalized Bessel functions $\mathcal{J}_k(x)$ is defined by
\begin{equation*}
\mathcal{J}_k(x):=\Gamma(k+1)\,(2/x)^k\,J_k(x).
\label{4} \end{equation*}
Here, $J_k(x)$ is the Bessel Function of the first kind \cite{Watson}.
Then
\begin{equation*}\label{}
 u(x)-\mathcal{M}_ru(x)=\int_{\R^d}\Big(1-\mathcal
{J}_{\gamma_k+d/2-1}(r|\xi| \Big)\;\mathcal{F}_ku(\xi)\;\mathcal{E}_k(i\xi,x)\; w_k(\xi)\;d\xi.
\end{equation*}
Substituting this expression in the integrand of \eqref{Intf}, one has
\begin{equation*}
(-\Delta_k)^{\alpha/2} u(x)=\frac{1}{\pi_{k,d}(\alpha)}\int_{0}^{\iy}\int_{\R^d}\frac{1-\mathcal
{J}_{\gamma_k+d/2-1}(r|\xi|)}{r^{1+\alpha}}\;\mathcal{F}_ku(\xi)\;\mathcal{E}_k(i\xi,x) \;w_k(\xi)\;d\xi.
\end{equation*}
Applying Fubini-Tonelli's theorem and taking into account of integral representation of the function $|\lambda|^{\gamma}$,  provided by \cite[Lemma 3.2]{B1},
\begin{equation*}\label{Bouza}
|\lambda|^\gamma=\frac{2^{\gamma+1}\Gamma(\nu+\tfrac{\gamma}{2}+1)}{\Gamma(\nu+1)
|\Gamma(-\tfrac{\gamma}{2})|}\int_{0}^{\infty}\big(1-\mathcal{J}_\nu(\lambda x)\big)\frac{dx}{x^{\gamma+1}},\quad 0<\gamma<2,
\end{equation*}
we obtain
\begin{align*}
(-\Delta_k)^{\alpha/2} u(x)&=\frac{1}{\pi_{k,d}(\alpha)}\int_{\R^d}\Big(\int_{0}^{\iy}\frac{1-\mathcal
{J}_{\gamma_k+d/2-1}(r|\xi|)}{r^{1+\alpha}}\,dr\Big)\mathcal{F}_ku(\xi)\mathcal{E}_k(i\xi,x) w_k(\xi)\,d\xi\\&=\frac{\Gamma(\gamma_k+\frac{d}{2})|\Gamma(-\frac{\alpha}{2})|}
{\pi_{k,d}(\alpha)2^{\alpha}\Gamma(\gamma_k+\frac{d+\alpha}{2})}
\int_{\R^d}|\xi|^\alpha\mathcal{F}_ku(\xi)\mathcal{E}_k(i\xi,x) w_k(\xi)\,d\xi.
\end{align*}
Therefore,
\begin{equation*}
\mathcal{F}_k (-\Delta_k)^{\alpha/2}u(\xi)=\frac{\Gamma(\gamma_k+\frac{d}{2})|\Gamma(-\frac{\alpha}{2})|}
{\pi_{k,d}(\alpha)2^{\alpha}\Gamma(\gamma_k+\frac{d+\alpha}{2})}|\xi|^\alpha \mathcal{F}_ku(\xi),\quad \xi \in \R^d.
\end{equation*}
In order to fulfill the equation \eqref{Fourf},  we impose that the normalized constant, in the above equation,  satisfies
\begin{equation*}
\frac{\Gamma(\gamma_k+\frac{d}{2})|\Gamma(-\frac{\alpha}{2})|}
{\pi_{k,d}(\alpha)2^{\alpha}\Gamma(\gamma_k+\frac{d+\alpha}{2})}=1.
\end{equation*}
For this to happen, we necessarily have
\begin{equation*}
\gamma_{\al,k,d}=\frac{c_k|\Gamma(-\frac{\al}{2})|}{2^{\alpha +\gamma_k+d/2}
\Gamma(\gamma_k+\tfrac{\alpha+d}{2})}.
\end{equation*}
\end{proofp}
\begin{lem} For $\alpha\in (0,\,2)$, we have
\begin{equation*}
\int_{\R^d}\frac{2-\mathcal{E}_k(i\xi,y)-\mathcal{E}_k(-i\xi,y)}{|y|^{\al+2\gamma_k+d}}w_k(y)dy= \frac{c_k|\Gamma(-\frac{\al}{2})||\xi|^{\alpha}}{2^{\alpha +\gamma_k+d/2} \Gamma(\gamma_k+\tfrac{\alpha+d}{2})}.
\end{equation*}
\end{lem}
\begin{proofp}
From the above property \eqref{f1} one has,
\begin{equation}\label{Hq3}
\int_{\R^d}\Gamma_k(t,y)\Big(2-\mathcal{E}_k(-i\xi,y)-\mathcal{E}_k(i\xi,y)\Big)\;w_k(y)\;
dy=2\;\Big(1-e^{-t|\xi|^2}\Big).
\end{equation}
Multiplying the members of \eqref{Hq3} by $t^{1+\al/2}$ and integrating over $(0,\,\iy)$ with respect to the variable $t$, we get
\begin{equation*}\label{Hq1}
\int_0^{\iy}\int_{\R^d}\frac{\Gamma_k(t,y)}{t^{1+\al/2}}\Big(2-\mathcal{E}_k(-i\xi,y)-\mathcal{E}_k(i\xi,y)\Big)
\;w_k(y)\;dy\;dt=2\int_{0}^{\iy}\frac{1-e^{-t|\xi|^2}}{t^{1+\al/2}}.
\end{equation*}
From the integral representation of the Dunkl Kernel given in \eqref{EV}, we deduce that
\begin{equation}\label{cos}
2-\mathcal{E}_k(iy,\xi)-\mathcal{E}_k(-iy,\xi)=2\int_{\R^d}\Big(1-\cos(\langle y,\eta\rangle)\Big)\mu^k_\xi(\eta)\geq 0.
\end{equation}
Taking into account \eqref{cos} and applying Fubini-Tonelli's theorem, we obtain
\begin{equation*}\label{Hq1}
\int_{\R^d}\int_{0}^{\iy}\frac{\Gamma_k(t,y)}{t^{1+\al/2}}\;dt\;
\Big(2-\mathcal{E}_k(-i\xi,y)-\mathcal{E}_k(i\xi,y)\Big)w_k(y)
\,dy\,=2\int_{0}^{\iy}\frac{1-e^{-t|\xi|^2}}{t^{1+\al/2}}dt.
\end{equation*}
A straightforward computation shows that
\begin{equation*}
\int_{0}^{\iy}\frac{\Gamma_k(t,y)}{t^{1+\frac{\alpha}{2}}}\,dt=\frac{\int_{0}^{\infty}\,
t^{-(1+\gamma+\alpha/2+d/2)}e^{-\frac{ |y|^2}{4t\, }}dt}{2^{\gamma_k+d/2}c_k}=2^{\alpha +\gamma_k+d/2}\frac{
\Gamma(\gamma_k+\tfrac{\alpha+d}{2})}{c_k|y|^{\alpha+2\gamma_k+d}}.
\end{equation*}
Therefore,
\begin{equation*}
\int_{\R^d}\frac{2-\mathcal{E}_k(i\xi,y)-\mathcal{E}_k(-i\xi,y)}
{|y|^{\al+2\gamma_k+d}}w_k(y)dy=\frac{c_k|\Gamma(-\frac{\al}{2})|}{2^{\alpha +\gamma_k+d/2}
\Gamma(\gamma_k+\tfrac{\alpha+d}{2})}.
\end{equation*}
\end{proofp}

\begin{ex}[The rank one case]  In the case $d=1$, the root system $\mathcal{R}$ equalt to $\{\pm \sqrt{2}\},$ $G=\mathbb{Z}_2$ and $w_k(x)=|x|^{2k}.$  Accordingly, the Dunkl operator
$T_k$, associated with the multiplicity parameter $k\geq  0$, is given by
\begin{equation*}
T_k :=\partial_x+\,\frac{k}{x}\,(1-s)\,,\quad \mbox{\it (see also \cite[Definition 4.4.2]{11}),}
\label{Dunkoul}
\end{equation*}
and the corresponding Dunkl Laplacian operator $\Delta_k$ is given by
\begin{align*}\label{Laplouce}
\Delta_k:=T_k^{2} =\partial^2_x +\frac{2k}{x}\partial_x -\frac{k}{x^2}(\,1\,-\,s\,),
\end{align*}
where $s$ is the reflection operator, which acts on a function $f(x)$ of real variable as:
\begin{equation*}
(sf)(x):=f(-x).
\end{equation*}
\end{ex}
Now, consider the so-called {\em nonsymmetric Bessel function}, also called {\em Dunkl-type Bessel function}, in the rank one case (see \cite[\S4]{11}):
\begin{equation*}
\mathcal{E}_k(x):=\mathcal{J}_{k-1/2}(ix)+\frac{\,x}{2k+1}\,\mathcal{J}_{k+1/2}(ix).
\label{5}
\end{equation*}
 Then, we have the eigenvalue equations
\begin{equation*}
T_k(\mathcal{E}_k(i\lambda x)) =i\,\lambda\,\mathcal{E}_k(i\lambda x),\quad \Delta_k(\mathcal{E}_k(i\lambda x))=-\lambda^2\, \mathcal{E}_k(i\lambda x),
\label{3} \end{equation*}
and the Dunkl transform is given by
\begin{equation*}
(\mathcal{F}_k f)(\lambda)=\widehat{f}(\lambda):=\frac{1}{2^{k+\frac{1}{2}} \Gamma(k+\frac{1}{2})} \int_{-\iy}^\iy f(x)\,\mathcal{E}_k(-i\lambda x)\, |x|^{2k}dx.
\end{equation*}
In \cite{Roesler1}, R\"{o}sler introduced  the following generalized translation $\tau^x$ defined by
\begin{eqnarray}\label{trs}
\tau^xu(y)&:=&\frac{1}{2}\int_{-1}^{1}u(\sqrt{x^2+y^2-2xyt})(1+\frac{x-y}{\sqrt{x^2+y^2-2xyt}})h_k(t)dt\\
&&\nonumber+\quad \frac{1}{2}\int_{-1}^{1}u(-\sqrt{x^2+y^2-2xyt})(1-\frac{x-y}{\sqrt{x^2+y^2-2xyt}}) \;h_k(t)\;dt,
\end{eqnarray}
where,
\begin{equation*}
h_k(t)=\frac{\Gamma(k+1/2)}{2^{2k}\sqrt{\pi}\Gamma(k)}(1+t)(1-t^2)^{k-1}.
\end{equation*}
\begin{teo} Let $u\in S(\R)$. Then we have
\begin{align*}
(-\Delta_k)^{\alpha/2}u(x) &=\frac{2^{\alpha} \Gamma(k+\tfrac{\alpha+1}{2})} {\Gamma(k+\frac{1}{2})
|\Gamma(-\tfrac{\alpha}{2})|}\;\lim_{\varepsilon\rightarrow0}\int_{|x|\geq \varepsilon}\int_{-1}^{1}
\tfrac{2u(x)-\varpi_{k,x}^{+}(t,y)u(\sqrt{x^2+y^2-2xyt})- \varpi_{k,x}^{-}(t,y)u(-\sqrt{x^2+y^2-2xyt})}{|y|^{\alpha+k+1/2}}dt dy,
\end{align*}
where \begin{equation*}
\varpi_{k,x}^{\pm}(t,y)=\frac{\Gamma(k+1/2)}{2^{2k}\sqrt{\pi}\Gamma(k)}
(1\pm\frac{x-y}{\sqrt{x^2+y^2-2xyt}})(1+t)(1-t^2)^{k-1}.
\end{equation*}
\end{teo}
\begin{proofp} In dimension one, the formula \eqref{Frr2} reads
\begin{equation*}
(-\Delta_k)^{\alpha/2}u(x)=\frac{2^{\alpha+1} \Gamma(k+\tfrac{\alpha+1}{2})}
{\Gamma(k+\frac{1}{2}) |\Gamma(-\tfrac{\alpha}{2})|}\;\lim_{\varepsilon\rightarrow0}\int_{|x|\geq \varepsilon}\int_{\mathbb{R}}\frac{u(x)-\tau^xu(y)}{|y|^{\alpha+k+1/2}} dy,
\end{equation*}
and the result follows after substituting the expression \eqref{trs} in the above formula.
\end{proofp}
\begin{ex}[Radial functions]
Recall that a function $u$ defined on $\mathbb{R}^d$ is  radial, if  there exists a function $u_0$ defined on $[0,\iy[$ such that $u(x)=u_0(|x|),\;\,  x\in \mathbb{R}^d.$ In polar coordinates $x=r\omega,$ the Dunkl Laplacian $\Delta_k$ is expressed as follows
\begin{equation}\label{rad}
\Delta_k=\frac{d^2}{dr^2}+\frac{2\gamma_k+d-1}{r}\frac{d}{dr}\Delta_{k,\mathbb{S}},
\end{equation}
where $\Delta_{k,\mathbb{S}}$ is the analogue of the Laplace-Beltrami operator on the sphere $\mathbb{S}^{d-1}.$ We refer the reader to \cite{Yu} for more details concerning  $\Delta_{k,\mathbb{S}}$.
From \eqref{rad}, we see that the operator $\Delta_k$ acts on the radial function $u(x)=u_0(|x|)$ as follows:
$$\Delta_u(x)=\mathcal{L}_{\gamma_k+d/2-1}u_0(|x|),$$
where the Bessel operator is
\begin{equation*}\label{Bessel1}
\mathcal{L}_{\gamma_k+d/2-1}=\frac{d^2}{dr^2}+\frac{2\gamma_k+d-1}{r}\frac{d}{dr}.
\end{equation*}
The Dunkl transform and the Fourier Bessel transform are deeply connected. In fact the Dunkl transform of a radial function $u\in L^1_k(\R^d)$  is again radial, i.e.
\begin{equation*}\label{radii}
\mathcal{F}_ku(x)=\mathcal{F}^B_{\gamma_k+d/2-1}u_0(|x|),
\end{equation*}
where $\mathcal{F}^B_{\gamma_k+d/2-1}u_0$ is  the Fourier-Bessel transform  of
$u_0\in L^1((0,\infty),d\sigma_\nu),$  given by
\begin{equation*}\label{fb1}
\mathcal{F}^B_{\gamma_k+d/2-1}u_0(r)= \frac{2^{1-\gamma_k-d/2}}{\Gamma(\gamma_k+d/2)}\int_{0}^\iy u_0(s)\,\mathcal{J}_{\gamma_k+d/2-1}(rs) \,\,s^{2\gamma_k+d-1}\,ds.
\end{equation*}
The fractional Bessel operator $(-\mathcal{L}_{\gamma_k+d/2-1})$ was considered in \cite{B1} and is defined as a Fourier-Bessel multiplier. Further, from \cite[Theorem 3.5]{B1}, we have
\begin{equation*}\label{R3}
(-\mathcal{L}_{\gamma_k+d/2-1})^{\alpha/2}f(r)= \frac{2^{\alpha+1}\Gamma(\gamma_k+\tfrac{d+\alpha}{2})}
{\Gamma(\gamma_k+\tfrac{d}{2}) |\Gamma(-\tfrac{\alpha}{2})|} \int_{0}^\iy
\frac{f(r)-T^rf(\varrho)}{\varrho^{1+\alpha}} \,d\varrho,
\end{equation*}
where
\begin{equation}\label{trr}
T^rf(\varrho)=\frac{\Gamma(\frac{d}{2})}{\sqrt{\pi}\Gamma(\frac{d-1}{2})}
\int_{0}^{\pi}f_0(\sqrt{r^2+\varrho^2\pm2r\varrho \cos \theta})\sin^{d-2} \theta d\theta
\end{equation}
\end{ex}
\begin{teo}Let $\al\in (0,2).$ For every radial function $u(x)=u_0(|x|) \in
S(\R^d)$, the fractional Laplace operator $(-\Delta_k)^{\alpha/2}u(x)$ is a radial function. Furthermore,
\begin{equation}\label{R3}
(-\Delta_k)^{\alpha/2}u(x)=\varsigma_k(\alpha)
\int_{0}^\iy\int_{0}^{\pi}\Big[u_0(|x|)-u_0(\sqrt{|x|^2+r^2\pm2r|x| \cos \theta}) \Big]
\frac{\sin^{2\gamma_k+d-1} \theta}{r^{1+\alpha}}\,d\theta dr,
\end{equation}
where the normalized constant $\varsigma_k(\alpha)$ is given by
\begin{equation}\label{xi}
\varsigma_k(\alpha)=\frac{2^{\alpha+1}\Gamma (\gamma_k+\frac{d+\alpha}{2})}{\sqrt{\pi}\Gamma
(\gamma_k+\frac{d}{2})|\Gamma (-\frac{\alpha}{2})|}.
\end{equation}
\end{teo}
\begin{proofp} From \eqref{radii}, the functions $\mathcal{F}_ku(\xi)$ and  $\mathcal{F}_k^{-1} \big(|\xi|^\alpha\mathcal{F}_ku(\xi)\big)(x)$ are radial. Then,
\begin{align*}
\mathcal{F}_k^{-1}\big(|\xi|^\alpha\mathcal{F}_ku(\xi)\big)(x)&=\mathcal{F}_k^{-1}\big(|\xi|^{\alpha}
\mathcal{F}^B_{\gamma_k+d/2-1}u_0(|\xi|) \big)(x) =\mathcal{F}^B_{\gamma_k+d/2-1}\big(r^{\alpha} \mathcal{F}^B_{\gamma_k+d/2-1}u_0(r) \big)(|x|),
\end{align*}
and
\begin{equation*}
(-\Delta_k )^{\alpha/2}u(x)=\frac{2^{\alpha+1}\Gamma(\gamma_k+\tfrac{d+\alpha}{2})}
{\Gamma(\gamma_k+\tfrac{d}{2}) |\Gamma(-\tfrac{\alpha}{2})|} \int_{0}^\iy
\frac{u_0(|x|)-T^{|x|}u_0(\varrho)}{\varrho^{1+\alpha}}\,d\varrho.
\end{equation*}
The result follows after substituting $T^{|x|}u_0(\varrho)$ in the above formula by its expression in \eqref{trr}.
\end{proofp}
\subsection{Bochner's subordination}
Our next is to derive a pointwise (integro-differential) formula for the fractional Dunkl-Laplacian opearator $(-\Delta_k)^{\alpha/2}$ by using the heat semigroup formalism. We begin with
a preliminary observation that connects the heat semigroup $\{e^{-t\Delta_k}\}_{t\geq 0}$ to the
spherical mean-value operator $\mathcal{M}_r^ku$.\par Recall that the heat semigroup $e^{-t\Delta_k}$ is given by \cite{Roesler2}
\begin{equation*}
\label{SH1}
e^{-t\Delta_k}f(x)=
\mathcal{F}_k^{-1}((e^{-t|\xi|^2}\mathcal{F}_kf(\xi))(x), \quad x\in\mathbb{R}^d.
\end{equation*}
Alternately \cite{Roesler2}
\begin{equation}\label{H1}
e^{-t\Delta_k}f(x)=\int_{\R^d}f(y)\tau^x\Gamma_k(t,y)w_k(y)dy.
\end{equation}

\begin{lem}
Let $u\in S(\R^d).$ For every $t\in (0,\,\iy),$ we have
\begin{equation*}\label{Es}
e^{-t\Delta_k}u(x)=u(x)+\frac{\sigma_k(d)}{(2t)^{\gamma_k+d/2}}\int_{0}^{\iy}e^{-r^2/4t}\Big[
\mathcal{M}_r^ku(x)-u(x)\Big]r^{2\gamma_k+d-1}dr.
\end{equation*}
\end{lem}
\begin{proofp} From \eqref{f1}, we have
\begin{eqnarray*}
e^{-t\Delta_k}u(x)-u(x)&=&\int_{\R^d}T^x\Gamma_k(y,t)\big[u(y)-u(x)\big]w_k(y)\,dy
=\int_{\R^d}\Gamma_k(y,t)\big[\tau^xu(y)-u(x)\big]w_k(y)\,dy\\
&=&\frac{1}{(2t)^{\gamma_k+d/2}}
\int_{0}^{\iy}e^{-r^2/4t} \int_{\mathbb{S}^{d-1}}\big[\tau^xu(r\omega )-u(x)\big]w_k(\omega)\,d\sigma(\omega)r^{2\gamma_k+d-1}dr\\
&=&\frac{\sigma_k(d)}{(2t)^{\gamma_k+d/2}}\int_{0}^{\iy}e^{-r^2/4t}\Big[
\mathcal{M}_r^ku(x)-u(x)\Big]r^{2\gamma_k+d-1}dr,
\end{eqnarray*}
which completes the proof.
\end{proofp}
\begin{prop}[Bochner representation]
For $0<\alpha<2$ and $u\in S(\R^d)$, the following holds
\begin{equation*}\label{Gsss}
(-\Delta_k)^{\alpha/2}u(x)=\frac{1}{|\Gamma(-\tfrac{\alpha}{2})|}
\int_{0}^{\iy} \big[e^{-t\Delta_k}u(x)-u(x)\big]\frac{dt}{t^{1+\tfrac{\alpha}{2}}}.
\end{equation*}\end{prop}
\begin{proofp} By Lemma 4.8, we have
\begin{eqnarray*}
\int_{0}^{\iy}\big[e^{-t\Delta_k}u(x)-u(x)\big]\frac{dt}{t^{1+\alpha/2}}&=&
\frac{\sigma_k(d)}{2^{\gamma_k+d/2}}\int_{0}^{\iy}\int_{0}^{\iy}
\frac{e^{-r^2/4t}}{t^{\gamma_k+\frac{d+\alpha}{2}+1}}\,dt\big[
\mathcal{M}_r^ku(x)-u(x)\big]r^{2\gamma_k+d-1}dr\,dt\\
&=&\frac{\sigma_k(d)}{2^{\gamma_k+d/2}}\int_{0}^{\iy}\Big(\int_{0}^{\iy}
\frac{e^{-r^2/4t}}{t^{\gamma_k+\frac{d+\alpha}{2}+1}}\,dt\Big)\big[
\mathcal{M}_r^ku(x)-u(x)\big]r^{2\gamma_k+d-1}dr,
\end{eqnarray*}
assuming that we can exchange the order of integration. Finally, since
\begin{equation}
\int_{0}^{\iy} \frac{e^{-r^2/4t}}{t^{\gamma_k+\frac{d+\alpha}{2}+1}}\,dt =2^{2\gamma_k+d+\alpha} \Gamma(\gamma_k+\frac{d+\alpha}{2})r^{-({2\gamma_k+d+\alpha)}},
\end{equation}
then,
\begin{equation*}
(-\Delta_k)^{\alpha/2}u(x)=\frac{1}{|\Gamma(-\tfrac{\alpha}{2})|}
\int_{0}^{\iy} \big[e^{-t\Delta_k}u(x)-u(x)\big]\frac{dt}{t^{1+\tfrac{\alpha}{2}}}.
\end{equation*}
\end{proofp}
\section{Extension problem}
\subsection{Fundamental solution of $(-\Delta_k)^{\alpha/2}$}
For $\alpha\in \mathbb{C},\;Re(\alpha)>0,$ we denote by $|x|^{-\alpha}$ the tempered distribution defined by
\begin{equation}\label{dis}
\langle |x|^{-\alpha},\,\varphi\rangle=\frac{1}{c_k} \int_{\mathbb{R}^d} \frac{\varphi(x)}{|x|^{\alpha}}\,w_k(x)dx, \quad \varphi \in \mathcal{S}(\mathbb{R}^d).
\end{equation}
For $\Re(\alpha)<d+2\gamma_k,$ the function $|x|^{-\alpha}$, as a distribution, generates a regular functional. However, if $\Re(\alpha)\geq d+2\gamma_k$, then $\;|x|^{-\alpha}$ is non-regular. In this case, it can be interpreted in the sense of regularization achieved by analytical continuation of the mapping $\alpha\rightarrow\langle|x|^{-\alpha},\,\varphi\rangle$ from the left half-plane $\Re(\alpha)<d+2\gamma_k$.
In the polar coordinates $y=r\omega$, the equation \eqref{dis} is expressed by a spherical mean type for the Dunkl operator in the following form
\begin{align*}
\langle |x|^{-\alpha},\,\varphi\rangle&=\frac{1}{c_k} \int_{0}^{\iy}\int_{\mathbb{S}^{d-1}}\varphi(ry)w_k(y)
\,d\sigma(y)\, r^{2\gamma_k+d-\alpha-1}\,dr\nonumber\\
&=\frac{1}{2^{\gamma_k+d/2-1}\Gamma(\gamma_k+d/2)}\int_{0}^{\iy} \mathcal{M}_{r}^k\varphi(x) r^{2\gamma_k+d-\alpha-1}dr, \quad \varphi \in \mathcal{S}(\mathbb{R}^d).\label{sphe}
\end{align*}
\begin{prop}
Let $\varphi\in \mathcal{S}(\mathbb{R}^d).$ The function $\alpha\rightarrow\langle |x|^{-\alpha},\,\varphi \rangle $ has an analytic extension to $\mathbb{C}-\{2\gamma_k+2p+d,\,\,p\in \mathbb{N}\}$, with simple poles at    $\alpha=2\gamma_k+2p+d$ and
\begin{equation}
\text{Res}(\langle|x|^{-\alpha},\,\varphi \rangle; 2\gamma_k+d+2p)= -\frac{\Delta_k^p\varphi(0)}{2^{2p+\gamma_k+d/2} \Gamma(\gamma_k+d/2+p)\,p!}.
\end{equation}
\end{prop}
\begin{proofp}  Observe that, for $\Re(\alpha)<2\gamma_k+2n+d$, we have
\begin{eqnarray}r
\int_{0}^{\iy} \mathcal{M}_{r}^k\varphi(0)r^{2\gamma_k+d-\alpha-1}dr&=&\int_{0}^{1} r^{2\gamma_k+d-\alpha-1}\Big(\mathcal{M}_{r}^k\varphi(0)-\sum_{p=0}^{n} \frac{r^{2p}\Delta_k^p\varphi(0)}{4^{p}(\gamma_k+d/2)_k\,p!}\Big)dr \label{t1}\nonumber\\
&&\quad +\sum_{p=0}^{n} \frac{\Delta_k^p\varphi(0)}{2^{2p}(\gamma_k+d/2)_p\,p!} \frac{1}{2\gamma_k+2p+d-\alpha}\nonumber  +\int_1^{\infty}\mathcal{M}_{r}^k\varphi(0)r^{2\gamma_k+d-\alpha-1}dr.
\end{eqnarray}
From Proposition \ref{mejme}  we have
\begin{equation*}\label{Remain}
\mathcal{M}_{r}^k\varphi(0)-\sum_{p=0}^{n} \frac{r^{2p}}{4^{p}(\gamma_k+d/2)_p\,p!}\Delta_k^p\varphi(0)=o(r^n),\qquad r\rightarrow 0.
\end{equation*}
The right hand-side of \eqref{t1} does not depend on the choice of $n$ ($n>\tfrac{1}{2} (\Re{(\alpha)}-2\gamma_k-d$) and from \eqref{Remain},  then the right hand-side of the formula \eqref{t1} yields an analytic continuation  of the mapping $\alpha\rightarrow\langle |x|^{-\alpha},\,\varphi \rangle $ on $\mathbb{C}-\{2\gamma_k+2p+d,\,\,p=0,1,\,2,\dots\}$, with simple poles at $\alpha=2\gamma_k+2p+d$ and \begin{equation*}
\text{Res}(\langle|x|^{-\alpha},\,\varphi \rangle;
2\gamma_k+2p+d)=-\frac{\sigma_k(d)\Delta_k^p\varphi(0)}{c_k2^{2p}(\gamma_k+d/2)_p\,p!}.
\end{equation*}
\end{proofp}
\begin{prop} The Dunkl transform of the distribution $|x|^{-\alpha}\in S^{'}({\mathbb{R}^d})$ is given by
\begin{equation*}
\widehat{|x|^{-\alpha}}=d_k(\alpha)
\left\{
     \begin{array}{l l}\label{c11}
      |x|^{\alpha-2\gamma_k-d},\quad \alpha \neq 2\gamma_k+d+2p,\,\,\alpha\neq-2p,\\ \\
       (-\Delta_k)^p \delta,\quad \quad \alpha=-2p.
       \end{array} \right.
\end{equation*}
Here $\delta$ is the Dirac delta function, $p\in \mathbb{N}\cup\{0\},$ the constant $d_k(\alpha)$ being given by
\begin{equation*}
d_k(\alpha)=
\left\{
        \begin{array}{l l}\label{c11}\frac{\Gamma(\gamma_k+(d-\alpha)/2)}
        {c_k 2^{\alpha-\gamma_k-d/2}\Gamma(\alpha/2)} ,\quad \alpha \neq 2\gamma_k+d+2p,\,\,\alpha\neq-2p,\\ \\
        1 ,\quad \quad \quad \quad  \quad \alpha=-2p.
        \end{array} \right.
\end{equation*}
\end{prop}
\begin{proofp}  Let $\varphi \in \mathcal{S}(\mathbb{R}^d),$ we have
\begin{equation*}\label{e1}
\langle \widehat{|x|^{-\alpha}},\,\varphi\rangle:=\langle|x|^{-\alpha},\, \widehat{\varphi}\rangle= \frac{1}{c_k} \int_{\mathbb{R}^d}\widehat{\varphi}(x)  |x|^{-\alpha}\,w_k(x)\,dx.
\end{equation*}
Using Parseval identity for the Dunkl transform, and taking account of \eqref{f1}, we obtain
\begin{equation}\label{p1}
\int_{\mathbb{R}^d} \mathcal{F}_k(\varphi)(x)e^{-t|x|^2} \,w_k(x)dx =\frac{1}{(2t)^{\gamma_k+d/2}} \int_{\R^d} \varphi(x) e^{-|x|^2/4t}\,w_k(x)dx.
\end{equation}
Now, multiplying both sides of the equation \eqref{p1} by $t^{-1+\alpha/2}$  and integrating over $(0,\,\infty)$ with respect to the variable $t$, we obtain, for $0<\alpha<2\gamma_k+d$, that
\begin{equation*}
\Gamma(\alpha/2)\int_{\R^d}\frac{ \widehat{\varphi}(x)}{|x|^{\alpha}} w_k(x)dx=
\frac{ \Gamma(\gamma_k+( d-\alpha)/2)}{2^{\alpha-\gamma_k-d/2}}\int_{\R^d} \frac{\varphi(x) }{ |x|^{2\gamma_k+d-\alpha}}\,w_k(x)dx.
\end{equation*}
By analytic continuation for $\alpha\in \mathbb{C}$ such that $\alpha \neq 2\gamma_k+d+2p,$ $p=0,\,1,\,2,\,\dots$ and $\alpha\neq -2p,\,\,p=0,\,1,\,2,\,\dots$, it is obvious that
\begin{equation}\label{zeta}
\widehat{|x|^{-\alpha}}=d_k(\alpha)|x|^{\alpha-2\gamma_k-d}, \quad
\mbox{where} \;d_k(\alpha)=\frac{\Gamma(\gamma_k+(d-\alpha)/2)}{c_k2^{\alpha-\gamma_k-d/2}\Gamma(\alpha/2)}.
\end{equation}
For  $\alpha=-2p$ and $p=0,\,1,\,\dots,$ we have
$$<\widehat{|x|^{2p}},\varphi>=<|x|^{2p},\widehat{\varphi}> =\frac{1}{c_k}\int_{\R^d}\widehat{\varphi}(x)|x|^{2p}w_k(x)dx =\frac{1}{c_k}\int_{\R^d}\widehat{(-\Delta_k)^{p}\varphi}(x)w_k(x)dx =\Delta_k ^{p}\varphi(0),$$
and finally,
$$\widehat{|x|^{2p}}=(-\Delta_k)^p \delta.$$
\end{proofp}

Recall $I_\nu$ the modified Bessel function of the first kind and $K_\nu$ the modified Bessel of the third kind, of order $\nu\neq 0,\,\pm1,\,\pm2,\,\dots,$ are given in \cite{Watson}, by
\begin{eqnarray*}
I_\nu(z)=\sum_{n=0}^{\iy}\frac{(z/2)^{\nu+2n}}{\Gamma(\nu+n+1)n!}\quad \mbox{and} \quad
K_\nu(z)=\frac{\pi}{2}\frac{I_{-\nu}(z)-I_{\nu}(z)}{\sin\nu z}, \quad |\arg(z)|<\pi.
\end{eqnarray*}
We will also consider
\begin{equation*}
q_{\alpha,y}(x)=\frac{d_k(\alpha)}{(y^2+|x|^2)^{\gamma_k+\frac{d-\alpha}{2}}}, \quad y>0,
\end{equation*}
and state the following preparatory result.
\begin{lem}\label{le} For $d\geq2$ and $y>0$, we have
\begin{enumerate}[1)]
\item $\widehat{q}_{\alpha,y}(x)=\frac{y^{\alpha/2}K_{\alpha/2}(y|x|)}{c_k2^{\alpha/2-1}
     \Gamma(\alpha/2)|x|^{\alpha/2}};$
\item $(-\Delta_k)^{\alpha/2}q_{\alpha,y}(x)=\frac{\Gamma(\gamma_k+(d+\alpha)/2)}{c_k2^{-\gamma_k-d/2}
      \Gamma(\alpha/2)}y^{\alpha}\big(y^{2}+|x|^2\big)^{-\gamma_k-\frac{d+\alpha}{2}}.$
\end{enumerate}
\end{lem}
\begin{proofp}  1) Let $\varphi\in S(\R^d)$, we have
\begin{equation} \label{a}
<\widehat{q}_{\alpha,y},\varphi> =d_k(\alpha)\int_{\R^d}\frac{\widehat{\varphi}(x)}{(y^2+|x|^2)^{\gamma_k+\frac{d-\alpha}{2}}}
w_k(x)dx=\frac{d_k(\alpha)}{\Gamma(\gamma_k+\frac{d-\alpha}{2})}
\int_0^{\iy}t^{\gamma_k+(d-\alpha)/2-1}e^{-ty^2}\int_{\R^d}
e^{-t|x|^2}\widehat{\varphi}(x)w_k(x)dx\,dt
\end{equation}
Using Parseval identity for the Dunkl transform and taking into account  \eqref{f1}, we obtain
\begin{equation}\label{a1}
\int_{\R^d}e^{-t|x|^2}\widehat{\varphi}(x)w_k(x)dx=\frac{1}{(2t)^{\gamma_k+d/2}}\int_{\R^d}
e^{-|x|^2/4t}\varphi(x)w_k(x)dx.
\end{equation}
Substituting \eqref{a1} in the integrand \eqref{a} and using Fubini-Tonelli's theorem, we obtain
\begin{equation*}
<\widehat{q}_{\alpha,y},\varphi>=\frac{d_k(\alpha)}{2^{\gamma_k+d/2}\Gamma(\gamma_k+\frac{d-\alpha}{2})}\int_{\R^d}\int_{0}^{\iy}
t^{-1-\alpha/2}e^{ty^2-|x|^2/4t}dt\,\varphi(x)w_k(x)dx
\end{equation*}
Now, using  formula \cite[p. 340]{Marichev}
\begin{equation*}
\int_{0}^{\iy}t^{\nu-1}e^{-(\frac{a}{t}+bt)}dt=2(a/b)^{\nu/2}K_\nu(2\sqrt{ab}),
\end{equation*}
we obtain
\begin{align*}
<\widehat{q}_{\alpha,y},\varphi>&=\frac{d_k(\alpha)y^{\alpha/2}}
{2^{\gamma_k+(d-\alpha)/2-1}\Gamma(\gamma_k+\frac{d-\alpha}{2})}
\int_{\R^d}|x|^{-\alpha/2}K_{\alpha/2}(y|x|)\,\varphi(x)w_k(x)dx
\\&=\frac{y^{\alpha/2}}{c_k2^{\alpha/2-1}\Gamma(\alpha/2)}
\int_{\R^d}|x|^{-\alpha/2}K_{\alpha/2}(y|x|)\,\varphi(x)w_k(x)dx,
\end{align*}
and this proves 1).\\

\noindent 2) From 1), we have
$$\widehat{(-\Delta_k)^{\alpha/2}q_{\alpha,y}}(\xi)=\frac{y^{\alpha/2}|\xi|^{\alpha/2}
K_{\alpha/2}(y|\xi|)}{c_k2^{\alpha/2-1}\Gamma(\alpha/2)}.$$
Applying inversion formula for the Dunkl transform,  we find that
\begin{equation*}
(-\Delta_k)^{\alpha/2}q_{\alpha,y}(x)=\frac{y^{\alpha/2}|x|^{1-\gamma_k-d/2}}
{c_k2^{\alpha/2-1}\Gamma(\alpha/2)}\int_{0}^{\iy}
K_{\alpha/2}(ry)J_{\gamma_k+d/2-1}(r|x|)\,r^{\gamma_k+(d+\alpha)/2}\,dr.
\end{equation*}
Making the substitutions $\nu\rightarrow \gamma_k+d/2-1,\,\,\mu\rightarrow \alpha/2$, $\lambda \rightarrow-\gamma_k-(d+\alpha)/2,\,\,a\rightarrow y,\,\,b \rightarrow |x|$ in the integral \cite[10.43.26]{11}
\begin{eqnarray*}
\int_{0}^{\infty}t^{-\lambda }K_{\mu}(at)J_{\nu }(bt)\,dt&=&
\frac{b^{\nu} \Gamma(\frac{1}{2}(\nu+\mu-\lambda+1))
\Gamma(\frac{1}{2}(\nu-\mu+\lambda+1))}{2^{\lambda+1 }a^{\nu-\lambda+1}\Gamma(\nu +1)}
 \hyp21{\frac{1}{2}(\nu+\mu-\lambda+1),\frac{1}{2}(\nu-\mu-\lambda+1)}
{\nu+1}{-\frac{a^{2}}{b^{2}}},
\end{eqnarray*}
where, $0<a<b,\,\, -1<\Re(\lambda)< \Re(\mu+\nu+1)$, we get
$$(-\Delta_k)^{\alpha/2}q_{\alpha,y}(x)= rac{\Gamma(\gamma_k +(d+\alpha)/2)}{c_k2^{-\gamma_k-d/2} \Gamma(\alpha/2)}y^{\alpha}\big(y^{2}+|x|^2 \big)^{-\gamma_k-\frac{d+\alpha}{2}}.$$
\end{proofp}
\begin{teo} Let $d\geq 2$, $0<\alpha<2,$ and recall the normalized constant $d_k(\alpha)$ is defined in \eqref{zeta}. Then, the function
\begin{equation*}
f_\alpha(x)=d_{k}(\alpha)|x|^{-(2\gamma_k+d-\alpha)},
\end{equation*}
is a fundamental solution for $(-\Delta_k)^{\alpha/2},$ that is $(-\Delta_k)^{\alpha/2}f_\alpha=\delta.$
\end{teo}
\begin{proofp}  Let $\varphi\in S(\R^d).$ Since $\lim_{y\rightarrow 0^+}q_{k,y}(x)=f_\alpha(x)$,  by Lebesgue dominated convergence theorem we have
\begin{equation*}
\lim_{y\rightarrow 0} \int_{\R^d}q_{k,y}(x)(-\Delta_k)^{\alpha/2}\varphi(x)w_k(x)dx =\int_{\R^d}f_\alpha(x)(-\Delta_k)^{\alpha/2}\varphi(x)w_k(x)dx.
\end{equation*}
On the other hand, form Lemma \eqref{le}, we can write
\begin{eqnarray*}
\int_{\R^d}(-\Delta_k)^{\alpha/2}q_{k,y}(x)\varphi(x)w_k(x)dx &=&\frac{\Gamma(\gamma_k
+(d+\alpha)/2)}{c_k2^{-\gamma_k-d/2}\Gamma(\alpha/2)}\int_{\R^d}\frac{y^{\alpha}\varphi(x)}{
(y^{2}+|x|^2 )^{\gamma_k+\frac{d+\alpha}{2}}}w_k(x)dx\\
&=&\frac{\Gamma(\gamma_k+(d+\alpha)/2)}{c_k2^{-\gamma_k-d/2}\Gamma(\alpha/2)}\int_{\R^d} \frac{\varphi(yx)}{
(1+|x|^2 )^{\gamma_k+\frac{d+\alpha}{2}}}w_k(x)dx.
\end{eqnarray*}
Using again Lebesgue dominated convergence theorem, we get
\begin{align*}
\lim_{y\rightarrow 0^+}\int_{\R^d}(-\Delta_k)^{\alpha/2}q_{k,y}(x)\varphi(x)w_k(x)dx
&=\varphi(0)\frac{\Gamma(\gamma_k +(d+\alpha)/2)}{c_k2^{-\gamma_k-d/2}\Gamma(\alpha/2)}\int_{\R^d}
\frac{1}{(1+|x|^2)^{\gamma_k+\frac{d+\alpha}{2}}}w_k(x)dx.
\end{align*}
To complete the proof of this Theorem,  it  suffices  to prove that
\begin{equation*}
\frac{\Gamma(\gamma_k +(d+\alpha)/2)}{c_k2^{-\gamma_k-d/2}\Gamma(\alpha/2)}\int_{\R^d}
\frac{1}{(1+|x|^2)^{\gamma_k+\frac{d+\alpha}{2}}}w_k(x)dx=1.
\end{equation*}
Indeed, using the polar coordinates $x=r\omega$, we get
\begin{equation*}
\int_{\R^d}\frac{1}{(1+|x|^2)^{\gamma_k+\frac{d+\alpha}{2}}}w_k(x)dx=
\sigma_k(d)\int_{0}^{\iy}\frac{r^{2\gamma+d-1}}{(1+r^2)^{\gamma_k+(d+\alpha)/2}}dr.
\end{equation*}
From the integral representation for the Beta function  (see, \cite[5.12.3]{11})
\begin{equation*}
\beta(a,b):=\frac{\Gamma(a)\Gamma(b)}{\Gamma(a+b)}=\int_{0}^{\iy}\frac{t^{a-1}}{(1+t)^{a+b}}dt,
\end{equation*}
we have
\begin{equation*}
\int_{\R^d}\frac{1}{(1+|x|^2)^{\gamma_k+\frac{d+\alpha}{2}}}w_k(x)dx=
\frac{c_k2^{-\gamma_k-d/2}\Gamma(\alpha/2)}{\Gamma(\gamma_k+(d+\alpha)/2)}.
\end{equation*}
\end{proofp}
\subsection{Extension problem}
The extension problem for the fractional Dunkl-Laplacian is a particular case of the general extension
problem proved in \cite{Sting1}, see also \cite{Sting2, M5}. In this section, we use the extension problem and the ideas of Caffarelli and Silvestre  \cite{Caffarelli}, and we show that the fractional Dunkl-Laplacian,  can be seen as the Dirichlet-Neumann map for a local degenerate elliptic equation. Let $u\in S(\R^d)$ and consider the following problem
\begin{equation}\label{pB}
\left\{
  \begin{array}{l l}\Delta_k U(x,y)+\frac{\partial^2U(x,y)}{\partial^2y}+\frac{1-\alpha}{y}\frac{\partial U(x,y)}{\partial y}=0,\quad (x,y)\in \R^{d}\times (0,\,\iy),\\
  \\
  U(x,0)=u(x),\quad x\in\R^d,\\
  \\
  U(x,y)\rightarrow 0, \quad y\rightarrow\iy.\end{array} \right.
\end{equation}
We define  the kernel $p_{\alpha,y}(x)$ on $\R^d$ by
\begin{equation*}\label{Poisson}
p_{\alpha,y}(x)=b_{k}(\alpha)\frac{y^\alpha}{\big(y^{2}+|x|^2\big)^{\gamma_k+\frac{d+\alpha}{2}}},\quad y>0,
\quad \mbox{\it where}\;\;
b_{k}(\alpha)=2^{\gamma_k+\tfrac{d}{2}}\frac{\Gamma(\gamma_k+\frac{d+\alpha}{2})}{c_k\Gamma(\tfrac\alpha2)}.
\end{equation*}
The following lemma will be useful in this section.
\begin{lem}
The kernel $p_{\alpha,y}$ has the following properties
\begin{enumerate}[1)]
\item $\displaystyle \widehat{p_{\alpha,y}}(\xi) = \frac{(y|\xi|)^{\alpha/2}}{c_k2^{\alpha/2-1}
\Gamma(\alpha/2)}K_{\alpha/2}(y|\xi|)$;
\item  $\displaystyle\|p_{\alpha,y}\|_{1,k}=\int_{\R^d}p_{\alpha,y}(x)w_k(x)dx=1.$
\end{enumerate}
\end{lem}
\begin{proofp}
The kernel $p_{\alpha,y}$ being a radial function, then by \eqref{radii}, one has
\begin{equation*}
\widehat{p_{\alpha,y}}(\xi)= \frac{b_{\alpha,k}y^\alpha}{|\xi|^{\gamma_k+d/2-1}}\int_{0}^\iy \frac{J_{\gamma_k+d/2-1}( r|\xi|)}{\big(y^{2}+r^2 \big)^{\gamma_k+\frac{d+\alpha}{2}}}\, \,\,r^{\gamma_k+d/2}\,dr.
\end{equation*}
The result follows from \cite[formula 6.565.4]{Marichev}.
\end{proofp}
\medskip

Note that  for $\alpha =1,$ the kernel \eqref{Poisson} takes the form
\begin{equation*}
p_{1,y}(x)=2^{\gamma_k+\tfrac{d}{2}}\frac{\Gamma(\gamma_k+\frac{d+1}{2})}{c_k\Gamma(\tfrac 12)}\frac{y}{\big(y^{2}+|x|^2 \big)^{\gamma_k+\frac{d+1}{2}}},\quad y>0.
\end{equation*}
which is the Poisson-Dunkl kernel for the half-space $\R^{d+1}_+,$ see \cite[formula (5.3)]{Roesler3}.
\begin{teo}
Let $u\in S(\R^d).$ Then, the solution $U$ to the extension problem \eqref{pB} is given by
\begin{equation*}\label{Poi}
U(x,y)=(p_{\alpha,y}\stk u)(x)=\frac{1}{c_k}\int_{\R^d}p_{\alpha,y}(\xi) \tau^{x}u(\xi)w_k(\xi)\,d\xi.
\end{equation*}
Furthermore, we have
\begin{equation}\label{f}
(-\Delta_k)^{\alpha/2}u(x)=-\frac{2^{\alpha-1}
\Gamma(\alpha/2)}{\Gamma(1-\alpha/2)}\lim_{y\rightarrow 0^+}y^{1-\alpha}\frac{\partial U}{\partial y}(x,y).
\end{equation}
\end{teo}
\begin{proofp}  Applying the Dunkl transform to the variable $x$ in \eqref{pB}, we find
\begin{equation}\label{pBB}
\left\{
  \begin{array}{l l}y^2Y''+(1-\alpha)yY'+|\xi|^2y^2Y=0,\\
  \\
  Y(0)=\mathcal{F}_ku(\xi),\quad \xi\in\R^d,\\
  \\
  Y(y)\rightarrow 0, \quad y\rightarrow0. \end{array} \right.
       \end{equation}
where $Y(y)=Y_\xi(y)=\mathcal{F}_kU(\xi,y).$ Thus, the general solution of \eqref{pBB} can be written in the form
\begin{equation*}
Y_\xi(y)=Ay^{\alpha/2} I_{\alpha/2}(y|\xi|)+By^{\alpha/2} K_{\alpha/2}(y|\xi|),
\end{equation*}
where $A$ and $B$ are constants depending on $\xi$. The condition $\lim_{\rightarrow 0}Y(y)$, gives $A =0$, and then
\begin{equation}\label{BB}
\mathcal{F}_kU(\xi,y)=By^{\alpha/2} K_{\alpha/2}(y|\xi|).
\end{equation}
To determine the constant $B$, we use the initial condition $\mathcal{F}_kU(\xi,0)=\mathcal{F}_ku(\xi),$ and the following asymptotic behavior of  $K_\nu(z)$ near zero:
\begin{equation*}
K_\nu(z)\sim 2^{\nu-1}\Gamma(\nu)z^{-\nu},
\end{equation*}
which yields
\begin{equation*}
\mathcal{F}_kU(\xi,y)\sim B2^{\alpha/2-1}\Gamma(\alpha/2)|\xi|^{-\alpha/2},\quad \text{as}\;y\rightarrow 0,
\end{equation*}
thus, we impose that
\begin{equation*}
B=\frac{|\xi|^{\alpha/2}\mathcal{F}_ku(\xi)}{2^{\alpha/2-1}\Gamma(\alpha/2)}.
\end{equation*}
Substituting this value of $B$ in \eqref{BB}, we  obtain
\begin{equation*}
\mathcal{F}_kU(\xi,y)=\frac{(y|\xi|)^{\alpha/2}}{2^{\alpha/2-1}\Gamma(\alpha/2)} K_{\alpha/2}(y|\xi|)\mathcal{F}_ku(\xi).
\end{equation*}
From Lemma 6.1, we get
\begin{equation*}
U(x,y)=p_\alpha\stk u(x),
\end{equation*}
therefore,
\begin{equation*}
U(x,y)=\frac{b_{\alpha,k}}{2c_k}\int_{\R^d}y^\alpha\frac{\tau^xu(\xi)+\tau^{-x} u(\xi)-2u(x)}
{\big(y^{2}+|\xi|^2 \big)^{\gamma_k+\frac{d+\alpha}{2}}}\; w_k(\xi)\;d\xi\;+ \;u(x).
\end{equation*}
Differentiating both sides of the above formula with respect to $y,$ we obtain
\begin{equation*}
y^{1-\alpha}\frac{\partial U}{\partial y}(x,y)= \frac{b_{\alpha,k}}{2c_k}\int_{\R^d} \frac{\alpha+(\gamma_k+d+\alpha)y^2}{\big(y^{2}+|\xi|^2\big)^{\gamma_k+\frac{d+\alpha}{2}}} \;\left[\tau^xu(\xi)+\tau^{-x} u(\xi)-2u(x)\right] \;w_k(\xi)\;d\xi,
\end{equation*}
and by Lebesgue dominated convergence theorem, we  find
\begin{align*}
\lim_{y\rightarrow 0^+}y^{1-\alpha}\frac{\partial U}{\partial y}(x,y)=
\frac{\alpha b_{\alpha,k}}{2c_k}\int_{\R^d}\frac{\tau^xu(\xi)+\tau^{-x}
u(\xi)-2u(x)} {|\xi|^{2\gamma_k+d+\alpha}}w_k(\xi)\,d\xi=\frac{\alpha b_{\alpha,k}}{c_k}\frac{\gamma_{k,d}(\alpha)}{2}(-\Delta_k)^{\alpha}.
\end{align*}
If, in the last equation, we  replace the expression \eqref{Cont} of the constant $\gamma_{k,d}(\alpha)$, we reach the conclusion that \eqref{f} is valid.
\end{proofp}

\end{document}